%
%
%
%

\hsize 159.2mm 
\vsize 246.2mm 
\font\Bbb=msbm10
  
\font\bigrm=cmr17  
\magnification=\magstep1
\def\C{\hbox{\Bbb C}}
\font\eightrm=cmr8
\font\eightbf=cmbx8
\font\eightit=cmti8  
\font\eightsl=cmsl8
\font\eightmus=cmmi8
\def\smalltype{\let\rm=\eightrm \let\bf=\eightbf
\let\it=\eightit \let\sl=\eightsl \let\mus=\eightmus
\baselineskip=9.5pt minus .75pt
\rm}

\centerline{\bigrm Continuity of the Complex Monge-Amp\`ere Operator}
\vskip .5in
\centerline{\sl Yang X\smalltype ING\footnote{$^\ast$}{Partially supported
by the Swedish
Natural Science Research Council.}} \vskip .8in 
\noindent{\bf 0. Introduction}
\bigskip\smallskip
\noindent Let $\Omega$ be an open subset in $\C^n$. $PSH(\Omega)$ will
stand for the set of all
plurisubharmonic (psh) functions on $\Omega$.
We use the standard notations $d=\partial +\overline\partial$ and
$d^c=i\,(\overline\partial-\partial)$.
The complex Monge-Amp\`ere operator $(dd^c)^n$ is,
via integrations by parts, well defined on
$PSH(\Omega)\,\cap\,L^\infty_{loc}(\Omega)$ and is 
continuous under monotone limits, that is, $(dd^cu_j)^n
\to (dd^cu)^n$ in the sense of currents if the monotone sequence of
functions $u_j$ 
converges to $u$ almost 
everywhere in $\Omega$, see [B-T2]. This basic fact implies an important
property that all psh
functons are quasi-continuous with respect to the capacity $C_n$ defined by 
$$C_n(E)=C_n(E,\Omega)=\sup\,\biggl\{\int\limits_E (dd^cu)^n;\, u\in
PSH(\Omega),\,0<u<1\biggr\} 
$$for any Borel set $E\subset \Omega,$ see [B-T2]. 

A natural question is to find the right notion of convergence $u_j\to u$
such that
$(dd^cu_j)^n \to (dd^cu)^n$ in the sense of currents. Cegrell [C] and
Lelong [L2] have observed that the convergence of $u_j$ to $u$
in $L^1_{loc}$ is not enough. In the papers [R] and [X], we used the
Hausdorff content, an outer measure
close to Hausdorff measure, to deal with this problem and obtained a
sufficient condition of 
the weak convergence $(dd^cu_j)^n \to (dd^cu)^n$. In section 1, by slightly
modifying the 
capacity $C_n$, we give a weaker condition to ensure  $(dd^cu_j)^n \to
(dd^cu)^n$. To see the 
sharpness of our conditions, we shall discuss properties of convergence of
functions $u_j$ to $u$ 
if the corresponding Monge-Amp\`ere measures converge in some sense.
Finally, in section 2 we give
an application of our results in the range of the Monge-Amp\`ere operator. 

The author would like to thank Urban Cegrell for helpful comments on an
earlier version of this paper.
\bigskip\bigskip\bigskip
\noindent{\bf 1. Continuity of the operator $(dd^c)^n$}
\bigskip\smallskip 
\noindent Similar to the paper [B-T2], we introduce an inner capacity
$C_{n-1}$ by
$$C_{n-1}(E)=C_{n-1}(E,\Omega)=\sup\,\bigl\{C_{n-1}(K);\, K\ is\ a\
compact\ subset\ of\ E\bigr\}$$
for each subset $E$ of $\Omega$, where for the compact subset $K$ we set
$$C_{n-1}(K)=\sup\,\biggl\{\int\limits_K (dd^cu)^{n-1}\wedge dd^c|z|^2;\,
u\in PSH(\Omega),\,0<u<1\biggr\}.$$

By the expansion of $\bigl(dd^c(u+|z|^2)\bigr)^n,$ we see that there exists
a constant $A_\Omega>0$ such that
$C_{n-1}(E)\leq A_\Omega C_n(E)$ for all subsets $E$ in $\Omega$. On the
other hand, using the
Chern-Levine-Nirenberg estimate, see [B-T2] and Lelong's result that 
$ r^{-2}\int_{\{|z-z_0|<r\}} (dd^cu)^{n-1}\wedge |z|^2$ for each $u\in
PSH(\Omega)$ and each 
$z_0$ in $\Omega$ is an 
increasing function of $r$, see [L1], we can easily prove
that for every $\Omega_1\subset\subset\Omega$ there exists a constant
$A_{\Omega_1}>0$ such that
$C_{n-1}(E)\leq A_{\Omega_1}\,H_2(E)$ holds for any subset $E$ of
$\Omega_1$, where $H_2$ denotes
Hausdorff 2-measure. This implies that $C_n$-capacity cannot be estimated
by $C_{n-1}$-capacity, since
there exists a compact                    
 subset $E$ in $\C^n$ with the Hausdorff dimension strictly less than 2 and
with the positive 
$C_n$-capacity (or equivalently, $E$ is not a pluripolar set).

Recall that a sequence of functions $u_j$ is said to converge to a function
$u$ in
$C_l$-capacity on a set $E$, where $l=n-1$ or $n$, if for each constant
$\delta>0$ we have 
$$C_l\bigl\{z\in E;\, |u_j(z)
-u(z)|>\delta\bigr\}\longrightarrow 0,\quad as\quad j\to \infty.$$
Using the quasi-continuity of psh functions with respect to $C_n$-capacity,
see [B-T2], we have
the following theorem.
 \bigskip
\noindent {\bf Theorem 1.}  $\,$\it Suppose that $u_j$ are locally
uniformly bounded psh functions 
in $\Omega$ and suppose that $u\in
PSH(\Omega)\,\cap\,L^\infty_{loc}(\Omega)$. Then the following 
assertions hold.
\smallskip
$(\romannumeral1)$\enskip If $u_j\to u$ in
$C_{n-1}$-capacity on 
each $E\subset\subset \Omega$, then $(dd^cu_j)^n\to (dd^cu)^n$ in the sense
of currents.
\smallskip
$(\romannumeral2)$\enskip If $u_j\to u$ in
$C_n$-capacity on 
each $E\subset\subset \Omega$, then $u_j(dd^cu_j)^n\to u(dd^cu)^n$ in the
sense of currents.
\rm 
\bigskip
\noindent{\it Proof.} We only give the proof of assertion
($\romannumeral1$), which in fact
also serves as the proof of assertion ($\romannumeral2$). 
We shall show by induction that for each positive integer $k\leq n$,
$(dd^cu_j)^k\to (dd^cu)^k$. It is clear for $k=1$ since the convergence
assumption implies that $u_j\to u$ in 
$L^1_{loc}(\Omega)$. Assume that it is true for $k=q<n.$ We 
prove now that $u_j(dd^cu_j)^q\to u(dd^cu)^q,$ which implies that the statement 
is true for $k=q+1$. For a given $\varepsilon>0$ we let $u=\phi+\psi$ on
$\Omega$ where 
$\phi$ is continuous, $\psi=0$ outside a subset of $\Omega$ with the
$C_n$-capacity
$<\varepsilon$, and the supremum norm of $\psi$ depends only on the
function $u$, see [B-T2]. Write   
$$u_j(dd^cu_j)^q-u(dd^cu)^q=(u_j-u)(dd^cu_j)^q+\psi
[(dd^cu_j)^q-(dd^cu)^q]+\phi [(dd^cu_j)^q-(dd^cu)^q].$$
The inductive assumption gives that the last term in the right hand side
converges to 0 in
the sense of currents. On the other hand, since
the positive measures $(dd^cu_j)^q\wedge (dd^c|z|^2)^{n-q}\leq \bigl(
dd^c(u_j+|z|^2)\bigr)^{n-1}
\wedge dd^c|z|^2$ which are majorized by $C_{n-1}$-capacity multiplied by a
constant 
not depending on $j$, it follows from the convergence assumption that the first
term in the right hand side also converges to 0 in the sense of currents.
Similarly, we can get that 
the second term makes arbitrarily small mass for all $j$ by choosing
$\varepsilon$ 
small enough. Therefore we have obtained the weak convergence
$u_j(dd^cu_j)^q\to u(dd^cu)^q,$ and the proof is complete.
\bigskip
Combining with Dini's theorem, Theorem 1 implies that the Monge-Amp\`ere
operator $(dd^c)^n$ is 
continuous under monotone limits. However, Theorem 1 is quite sharp as the
following result shows.
\bigskip
\noindent {\bf Theorem 2.} \it Let $E\subset\subset\Omega\subset \C^n$.
Suppose that $u_j$ are
locally uniformly bounded psh functions in $\Omega$ and suppose that there
exists a function $u\in PSH(\Omega)\,\cap\,L^\infty_{loc}
(\Omega)$ such that $u_j=u$ on $\Omega\setminus E$ for all $j$. Then the
following assertions
hold.
\smallskip
$(\romannumeral1)$\enskip $u_j\to u$ in $C_n$-capacity on $\Omega$ if and
only if
$u_j(dd^cu_j)^n,\ u(dd^cu_j)^n$ and $u_j(dd^cu)^n$ converge to $u(dd^cu)^n$
in the sense of currents.
\smallskip
$(\romannumeral2)$\enskip $u_j\to u$ in $C_{n-1}$-capacity on $\Omega$ if
and only if
$(dd^cu_j)^n,\ (dd^cu_j)^{n-1}\wedge dd^cu$ and $(dd^cu)^{n-1}\wedge
dd^cu_j$ converge to 
$(dd^cu)^n$ in the sense of currents.
\smallskip
$(\romannumeral3)$\enskip In the special case $n=2$, we have that $u_j\to
u$ in $C_1$-capacity 
on $\Omega$ if and only if $(dd^cu_j)^2\to (dd^cu)^2$ in the sense of currents
and $u_j\to u$ in $L^1_{loc}(\Omega)$.
\smallskip
$(\romannumeral4)$\enskip Under the additional assumption that either
$u_j\geq u$ or $u_j\leq u$ holds
in $\Omega$ for each $j$,
we have that $u_j\to u$ in $C_{n-1}$-capacity on $\Omega$ if and only if
$(dd^cu_j)^n\to (dd^cu)^n$ 
in the sense of currents.
\rm \bigskip
\noindent{\it Proof.} All of the ``only if$\,$" parts follow from the proof
of Theorem 1.
We only need to show the ``if$\,$" parts. Assume that $\delta$ is a
positive constant and assume
that $w$ is a function in $PSH(\Omega)$ with $0<w<1$. Choose an open set $\Omega^\prime$ such that
$E\subset\subset \Omega^\prime\subset\subset \Omega$.

$(\romannumeral1)$\enskip Using an integration by parts
and the Schwarz inequality, we have
$$\int\limits_{\{|u_j-u|>\delta\}} (dd^cw)^n\leq {1\over \delta^2}\,
\int\limits_{\Omega^\prime} (u_j-u)^2 (dd^cw)^n= -{1\over \delta^2}\, 
\int\limits_{\Omega^\prime} d(u_j-u)^2 \wedge d^cw\wedge (dd^cw)^{n-1}$$
$$   
\leq A_1\, \biggl(\int\limits_{\Omega^\prime} d(u_j-u)^2 \wedge
d^c(u_j-u)^2\wedge (dd^cw)^{n-1}\biggr)^{1/2}
$$ $$\leq  2A_1A_2\,
\biggl(\int\limits_{\Omega^\prime}d(u_j-u) 
\wedge d^c(u_j-u)\wedge (dd^cw)^{n-1}\biggr)^{1/2}, $$
where the
constant $A_1=1/\delta^2\,\bigl(\int_{\Omega^\prime} dw\wedge d^cw\wedge
(dd^cw)^{n-1}\bigr)^{1/2}$ is
uniformly bounded for all functions $w\in PSH(\Omega)$ with $0<w<1$,
see the Chern-Levine-Nirenberg estimate in [B-T2], and the constant $A_2$
exceeds $|u_j(z)-u(z)|$ for all $j$ and $z\in\Omega$.
Again by an integration by parts, the last integral in the above inequality
is equal to
$$\int\limits_{\Omega^\prime}d(u_j-u)\wedge d^cw\wedge dd^cu_j\wedge
(dd^cw)^{n-2}-
\int\limits_{\Omega^\prime}d(u_j-u)\wedge d^cw\wedge dd^cu\wedge
(dd^cw)^{n-2}.$$
Applying the Schwarz inequality in each term of this difference, it then
turns out from the trivial
inequalities $dd^cu_j\leq dd^c(u_j+u)$ and $dd^cu\leq dd^c(u_j+u)$ that the
last difference does not
exceed
$$A_3\,\biggl( \int\limits_{\Omega^\prime} d(u_j-u)\wedge d^c(u_j-u)\wedge
dd^c(u_j+u)\wedge (dd^cw)^{n-2}\biggr)^{1/2},$$
where the constant $A_3$ does not depend on $w$ and $j$ because of the
Chern-Levine-Nirenberg 
estimate. Now we get an integral involving $(dd^cw)^{n-2}$. ( Observe that
we began with the integral
involving $(dd^cw)^{n-1}$. )  
We repeat this argument $n-2$ more times and finally find a constant $A_4$,
not depending on $w$ and $j$, such that 
$$\int\limits_{\{|u_j-u|>\delta\}} (dd^cw)^n\leq A_4\,\biggl(
\int\limits_{\Omega^\prime} d(u_j-u)
\wedge d^c(u_j-u)\wedge \bigl(dd^c(u_j+u)\bigr)^{n-1}\biggr)^{1/2^n}$$
$$\leq A_4\, (n!)^{1/2^{n-1}}\,\biggl( \int\limits_{\Omega^\prime} d(u_j-u)
\wedge d^c(u_j-u)\wedge \sum\limits_{k=0}^{n-1}(dd^cu_j)^{n-1-k}\wedge
(dd^cu)^k \biggr)^{1/2^n}$$
$$= A_4\, (n!)^{1/2^{n-1}}\,\biggl( \int\limits_{\Omega^\prime} (u-u_j)
 (dd^cu_j-dd^cu)\wedge \sum\limits_{k=0}^{n-1}(dd^cu_j)^{n-1-k}\wedge
(dd^cu)^k \biggr)^{1/2^n}$$
$$= A_4\, (n!)^{1/2^{n-1}}\,\biggl( \int\limits_{\Omega^\prime} (u-u_j)
\bigl( (dd^cu_j)^n- (dd^cu)^n\bigr) \biggr)^{1/2^n}.$$
Since the function $u-u_j$ has a compact support in $\Omega^\prime$, it
follows from
the convergence assumptions that the last integral converges to $0$ as
$j\to\infty$. Hence we have proved 
that for any $\delta>0$
$$\lim\limits_{j\to\infty} \biggl[ \ \sup\ \biggl\{ \ \int\limits_{\{
|u_j-u|>\delta \}}(dd^cw)^n;\, 
w\in PSH(\Omega),\, 0<w<1
\biggr\}\ \biggr]=0,$$
which completes the proof of assertion ($\romannumeral1$). 

$(\romannumeral2)$\enskip If we begin with the integral
$\int_{\{|u_j-u|>\delta\}}
(dd^cw)^{n-1}\wedge dd^c|z|^2$ and use the same argument as in the proof of
$(\romannumeral1)$, we
can get the ``if$\,$" part of assertion $(\romannumeral2)$.

$(\romannumeral3)$\enskip By the quasi-continuity of the function $u$ 
with respect to $C_2$-capacity and the convergence
assumption that $u_j\to u$ in $L^1_{loc}(\Omega)$, we can easily get the
weak convergence
 $dd^cu_j\wedge dd^cu\to (dd^cu)^2$.
Hence assertion $(\romannumeral3)$ follows directly from assertion
$(\romannumeral2)$.

$(\romannumeral4)$\enskip Similar to the proof of $(\romannumeral1)$ , we
can find a constant $A_5$,
not depending on $w$ and $j$, such that
$$\int\limits_{\{|u_j-u|>\delta\}} (dd^cw)^{n-1}\wedge dd^c|z|^2$$  $$\leq
 A_5\,\biggl( \int\limits_{\Omega^\prime} d(u_j-u)
\wedge d^c(u_j-u)\wedge \bigl(dd^c(u_j+u)\bigr)^{n-2}\wedge
dd^c|z|^2\biggr)^{1/2^{n-1}}$$
$$=A_5\,\biggl( \int\limits_{\Omega^\prime} (u-u_j)
dd^c(u_j-u)\wedge \bigl(dd^c(u_j+u)\bigr)^{n-2}\wedge
dd^c|z|^2\biggr)^{1/2^{n-1}}.$$
Since $u-u_j$ does not change sign on $\Omega$ for every $j$, the last
integral is majorized by
$$(n!)^2\,\bigg| \int\limits_{\Omega^\prime} (u_j-u)dd^c|z|^2\wedge
\sum\limits_{k=0}^{n-1}(dd^cu_j)^{n-1-k}\wedge (dd^cu)^k
\bigg|=(n!)^2\,\bigg|\int\limits_{\Omega^\prime}
|z|^2\bigl((dd^cu_j)^n-(dd^cu)^n\bigr)\bigg|
,$$ which, by the assumptions, converges to $0$ as $j\to\infty$ and hence 
assertion $(\romannumeral2)$ has been proved. We have thus completed the
proof of Theorem 2.
\bigskip
For the ``if$\,$" parts of Theorem 2, we require that all functions $u_j$
coincide with $u$ 
outside a relatively compact subset of $\Omega$.
This requirement cannot be replaced  by the weaker restriction that all $u_j$
have the same boundary values on $\partial\Omega$ as the function $u$. For
instance, psh functions
$u_j(z)=\max\,(j\ln |z|,\, -1)$ converge to the function $u(z)\equiv 0$
nowhere inside the open
unit ball in $\C^n$ as $j\to\infty$.
However, every function $u_j$ for $j>0$ vanishes on the boundary of the
unit ball, and its Monge-Amp\`ere measure $(dd^cu_j)^n$ is a constant
multiple of Lebesgue
measure on the sphere $|z|=e^{-1/j}$, which implies that together with the
$u$, the functions $u_j$
satisfy the other assumptions of the ``if$\,$" parts of all assertions in
Theorem 2, except assertion 
$(\romannumeral3)$. In the following we shall give a slightly weaker
condition on functions
near the boundary instead of that given in Theorem 2. To do this we prefer
to set up the following
inequality, which we feel has some interest in itself.
\bigskip
\noindent {\bf Lemma 1. } \it Let $\Omega$ be a bounded open subset in
$\C^n$ and let $u,\,v\in
PSH(\Omega)\,\cap\,L^\infty(\Omega)$ satisfy $\lim\inf_{z\to\partial\Omega}
\bigl(
u(z)-v(z)\bigr)\geq 0.$ Then for any constant $r\geq 1$ and all $w_j\in
PSH(\Omega)$ with $0\leq w_j
\leq 1,\ j=1,2,\dots, n$, we have $${1\over {(n!)^2}} \int\limits_{\{u<v\}}
(v-u)^n\
dd^cw_1\wedge\cdots\wedge dd^cw_n+\int\limits_{\{u<v\}} (r-w_1)\
(dd^cv)^n\leq \int\limits_{\{u<v\}}
(r-w_1)\ (dd^cu)^n. $$

\rm \bigskip
\noindent{\it Proof.} We first prove Lemma 1 for continuous functions $u$
and $v$ in
$\Omega$. In this case we can assume, without loss of generality, that
$\Omega=\{u<v\}$. For 
each constant $\varepsilon>0$ we define a function $v_{\varepsilon }=\max
(u,\ v-\varepsilon)$, which 
converges increasingly to $v$ in $\Omega$ as $\varepsilon\searrow 0$ and
coincides with $u$ near the 
boundary $\partial\Omega$. So integrations by parts yield 
$$ \int\limits_{\Omega} (v_\varepsilon-u)^n\ dd^cw_1\wedge\cdots\wedge
dd^cw_n = 
\int\limits_{\Omega}(w_n-1)\ dd^c\bigl((v_\varepsilon-u)^n\bigr)\wedge
dd^cw_1\wedge\cdots\wedge dd^cw_{n-1} $$ 
$$=n(n-1)\int\limits_{\Omega}(w_n-1)\,(v_\varepsilon-u)^{n-2}\
d(v_\varepsilon-u)\wedge  d^c(v_\varepsilon-u)\wedge
dd^cw_1\wedge\cdots\wedge dd^cw_{n-1}$$
$$+n\int\limits_{\Omega}(w_n-1)
(v_\varepsilon-u)^{n-1}\ dd^c(v_\varepsilon-u)\wedge
dd^cw_1\wedge\cdots\wedge dd^cw_{n-1}$$ 
$$\leq
n\int\limits_{\Omega}(v_\varepsilon-u)^{n-1}\ dd^c(v_\varepsilon+u)\wedge
dd^cw_1\wedge\cdots\wedge dd^cw_{n-1}.$$
Repeating this process $n-2$ more times, we have 
$$ \int\limits_{\Omega} (v_\varepsilon-u)^n\ dd^cw_1\wedge\cdots\wedge
dd^cw_n\leq n! 
\int\limits_{\Omega}(v_\varepsilon-u)\
\bigl(dd^c(v_\varepsilon+u)\bigr)^{n-1}\wedge dd^cw_1$$
$$ \leq (n!)^2 
\int\limits_{\Omega}(v_\varepsilon-u)\
dd^cw_1\wedge\sum\limits_{k=0}^{n-1}\bigl(dd^cv_\varepsilon\bigr)^{n-1-k}
\wedge \bigl(dd^c
u\bigr)^k $$
$$=(n!)^2
\int\limits_{\Omega}(r-w_1)\
\bigr(dd^cu\bigr)^n-(n!)^2\int\limits_{\Omega}(r-w_1)\
\bigr(dd^cv_\varepsilon\bigr)^n. $$
But $(r-w_1)\,(dd^cv_\varepsilon)^n\longrightarrow (r-w_1)\,(dd^cv)^n$ as
currents when
$\varepsilon\searrow 0$, see [B-T2] and hence we have obtained the required
inequality for
continuous functions $u$ and $v$.

The general case will then follow by an approximation argument. As in the
proof of Theorem 4.1 in [B-T2],
we may assume that there exists an open $E\subset\subset\Omega$ such that
$u(z)-v(z)\geq \delta>0$ for all
$z\in \Omega\setminus E$. Otherwise, replace $u$ by $u+2\delta$ and then
let $\delta\searrow
0$. We can thus choose two decreasing sequences of smooth psh functions $u_k$ and $ v_j$ in a neighbourhood $\Omega^\prime$
of $\overline E$ such that $\lim_{k\to\infty} u_k=u,\ \lim_{j\to\infty}
v_j=v$ in $\Omega^\prime$ and $u_k\geq v_j$ near
the boundary $\partial\Omega^\prime$. For smooth functions $u_k$ and $v_j$
in $\Omega^\prime$, we have proved
the following inequality
$${1\over {(n!)^2}} \int\limits_{\Omega^\prime} \chi_{{\{z\in
\Omega^\prime;\, u_k<v_j\}}} (v_j-u_k)^n\
dd^cw_1\wedge\cdots\wedge dd^cw_n+\int\limits_{\{z\in \Omega^\prime;\,
u_k<v_j\}} (r-w_1)\ (dd^cv_j)^n$$  $$
\leq \int\limits_{\{z\in \Omega^\prime; \, u_k <v_j\}} (r-w_1)\ (dd^cu_k)^n,$$
where $\chi_{_E}$ denotes the characteristic function of a set $E$. Letting
$j\to\infty$ and 
then $k\to\infty$ and using Fatou Lemma, we get that the limit inferior of
the first term in the left hand 
side exceeds 
$${1\over {(n!)^2}} \int\limits_{\{u<v\}}  (v-u)^n\
dd^cw_1\wedge\cdots\wedge dd^cw_n.$$
To handle the other two terms in the same inequality when $j\to\infty$ and
$k\to\infty$, we first observe that
$ (r-w_1)\ (dd^cu_k)^n \longrightarrow (r-w_1)\ (dd^cu)^n $ and $(r-w_1)\
(dd^cv_j)^n\longrightarrow 
(r-w_1)\ (dd^cv)^n$ as currents, see [B-T2]. Completely repeating the proof
of Theorem 4.1 in
[B-T2], one can get 
$${1\over {(n!)^2}} \int\limits_{\{u<v\}} (v-u)^n\
dd^cw_1\wedge\cdots\wedge dd^cw_n+\int\limits_{\{u<v\}} (r-w_1)\
(dd^cv)^n\leq \int\limits_{\{u\leq v\}}
(r-w_1)\ (dd^cu)^n. $$
Finally, applying the last inequality for functions $u+\delta_1$ instead of
$u$ and letting $\delta_1\searrow
0$, we get the required inequality and hence the proof is complete. 
\bigskip
\noindent {\it Remark.\rm} If the both sides of the inequality in Lemma 1
are devided by the
constant $r$ and then letting $r\to\infty$, we obtain the inequality 
$$\int\limits_{\{u<v\}} (dd^cv)^n\leq \int\limits_{\{u<v\}} (dd^cu)^n, $$
which is the result of the comparison theorem for the complex
Monge-Amp\`ere operator, due to Bedford and Taylor
[B-T2]. Our inequality also implies the following useful estimate.
\bigskip
\noindent {\bf Lemma 2. } \it Let $\Omega$ be a bounded open set and let
$u,\,v\in PSH(\Omega)\,\cap\,L^\infty(\Omega)$ satisfy $\
\lim\sup_{z\to\partial\Omega} | u(z)-v(z)|= 0$. Then the following inequality 
$$C_n\bigl\{ |u-v|\geq\delta \bigr\}\leq {{(n!)^2}\over {(1-k)^n\delta^n}}\
\bigl|\bigl|
(dd^cu)^n-(dd^cv)^n\bigr|\bigl|_{\{ |u-v|>k\delta \}}$$
holds for all constants $\delta>0$ and $0<k<1$, where 
$||\mu ||_E$ denotes the mass on $E$ of the total variation of a signed measure
$\mu$.
\rm
\bigskip
\noindent{\it Proof.} Since $|u-v\pm k\delta|\geq (1-k)\,\delta$ on the set
$\{|u-v|\geq\delta\}$,
we deduce from Lemma 1 that for $w_1=w_2=\dots=w_n=w\in PSH(\Omega)$ with
$\ 0<w<1$ 
$$\int\limits_{\bigl\{ |u-v|\geq\delta \bigr\}}(dd^cw)^n
$$  $$\leq {1\over {(1-k)^n\delta^n}}\ \Biggl[\ \int\limits_{\{u+\delta
\leq v\}} (v-u-k\delta)^n\,(dd^cw)^n
+ \int\limits_{\{v+\delta \leq u\}} (u-v-k\delta)^n\,(dd^cw)^n\Biggr]$$
$$\leq {1\over {(1-k)^n\delta^n}}\ \Biggl[\ \int\limits_{\{u+k\delta <v\}}
(v-u-k\delta)^n\,(dd^cw)^n
+ \int\limits_{\{v+k\delta <u\}} (u-v-k\delta)^n\,(dd^cw)^n\Biggr]$$
$$\leq {(n!)^2\over {(1-k)^n\delta^n}}\ \Biggl[\
\int\limits_{\{|u-v|>k\delta\}} (1-w)\,\bigl(
\chi_{\{u+k\delta <v\}}- \chi_{\{v+k\delta <u\}}
\bigr)\,\Bigl((dd^cu)^n-(dd^cv)^n\Bigr)\Biggr]$$
$$\leq {(n!)^2\over {(1-k)^n\delta^n}}\ \bigl|\bigl|
(dd^cu)^n-(dd^cv)^n\bigr|\bigl|_{\{|u-v|>k\delta\}},$$
which completes the proof.
\bigskip
As a direct consequence of Lemma 2 we have
\bigskip
\noindent {\bf Theorem 3. } \it Suppose that $\Omega$ is a bounded open set
and suppose that $u_j,\,u\in PSH(\Omega)\,\cap\,L^\infty(\Omega)$.
If $$\leqno{\qquad\quad (\romannumeral1)}\qquad\quad
\limsup\limits_{z\to\partial\Omega}\,| u_j(z)-u(z)|= 0 \quad {uniformly\ in\ j,}
\ {and} $$ $$\leqno{\qquad\quad (\romannumeral2)}\qquad\quad
\big|\big|(dd^cu_j)^n-(dd^cu)^n\big|\big|_E\longrightarrow 0\quad{ for\
any\ subset\ E\subset
\subset\Omega},$$
then $u_j\to u$ in $C_n$-capacity on $\Omega$.
\rm 
\bigskip 
Note that the uniformly vanishing condition ($\romannumeral1$) of Theorem 3
may be replaced by  $\limsup_{z\to\partial\Omega}\,| u_j(z)-u(z)|= 0$ for
each $j$, if we assume 
$||(dd^cu_j)^n-(dd^cu)^n||_{\Omega}\longrightarrow 0$ instead of condition
($\romannumeral2$).
Otherwise, we cannot weaken condition ($\romannumeral1$) in such a way, as
can be seen from the simple example
$u_j(z)=\max\,(j\ln |z|,\, -1)$ and $u(z)\equiv 0$ in the unit ball. It
should be also mentioned that 
condition ($\romannumeral2$) of Theorem 3 cannot be replaced by the weak
convergence $(dd^cu_j)^n\to (dd^cu)^n$
either, as the following example shows.
\bigskip
\noindent{\it Example.} Assume that $\Omega$ is an open unit ball in
$\C^1$. By Lemma 2 in [C] 
there exists a sequence of subharmonic functions $f_j$ with $-1\leq f_j\leq
1/2$ in $\Omega$ such that $f_j$
converges to a subharmonic function $f$ in the topology of
$L^1_{loc}(\Omega)$, but $f_jdd^cf_j$
does not converge to $fdd^cf$ as currents. Hence it follows from the
monotone convergence theorem in
[B-T2] that if the constant $A$ is big enough, $u_jdd^cu_j$ does not
converge to $udd^cu$ as currents,
where the functions $u_j=\max( A\ln |z|,\, f_j)$ and $u=\max( A\ln |z|,\,
f)$ coincide 
outside a compact subset of $\Omega$. By Theorem 1 we then have that $u_j$
does not converge 
to $u$ in $C_1$-capacity on
$\Omega$. However, $u_j\to u$ in $L_{loc}^1(\Omega)$ which gives the weak
convergence
$dd^cu_j\to dd^cu$.
\bigskip
In fact, we have an analogue of Theorem 3 for the convergence in
$C_{n-1}$-capacity. But we omit
the details here, because both the formulation and the proof are completely
similar.
\bigskip\bigskip\bigskip
\noindent{\bf 2. Range of the operator $(dd^c)^n$ }
\bigskip\smallskip
\noindent In this section we will give
an application of our results in the range of the Monge-Amp\`ere operator.
We shall study the problem: Find a function $u\in PSH(\Omega)\cap
L^\infty(\Omega)$ such that 
$(dd^cu)^n=\mu$ on $\Omega$, where $\mu$
is a given positive measure in the bounded open set $\Omega$. One necessary
condition of existence of such a solution is that 
there exists a subsolution for this problem, that is, there exists $v\in
PSH(\Omega)\cap L^\infty(\Omega)$
such that $(dd^cv)^n\geq \mu$.  In [C-S]
it was shown that for $\mu=fd\lambda$, where $f\in L^1(\Omega)$ and
$d\lambda$ is the Lebesgue measure,
the problem has a solution if there exists a subsolution. Now our result is
the following.
\bigskip
\noindent {\bf Theorem 4. }\it Suppose that there exists a function $v\in
PSH(\Omega)\cap 
L^\infty(\Omega)$ such that $(dd^cv)^n\geq \mu$ and suppose that there
exist Monge-Amp\`ere measures 
$\mu_j=(dd^cu_j)^n$ in $\Omega$ such that $||\mu_j-\mu||_\Omega\to 0$ as
$j\to\infty$, where
all functions $u_j\in  PSH(\Omega)\cap C(\,{\overline \Omega}\,)$ take the
same boundary values
on $\partial\Omega$. Then there exists a function $u\in PSH(\Omega)\cap
L^\infty(\Omega)$ such
that $(dd^cu)^n=\mu$ in $\Omega$.\rm
\bigskip
\noindent{\it Proof.} By passing to a subsequence, we may assume that  
$$2^n(n!)^2||\mu_j-\mu||_{\Omega}\leq {1\over 2^{(n+2)j}}\qquad for\ 
j=1,2,\dots.$$
So Lemma 2 gives that for any $\delta
>0$ $$C_n\bigl\{ |u_{j+1}-u_j|\geq \delta\bigr\}\leq {{2^n(n!)^2}\over
>{\delta^n}}\, 
||\mu_{j+1}-\mu_j||_{\Omega}$$
$$\leq {{2^n(n!)^2}\over {\delta^n}}\,\bigl( ||\mu_{j+1}-\mu||_{\Omega}+
||\mu-\mu_j||_{\Omega}\bigr)\leq {1\over {\delta^n 2^{(n+1)j}}}\qquad for\ 
j=1,2,\dots.$$
Choose a constant $A$ such that $A\geq |z|$ for all $z\in\Omega$ and choose
a constant $c$ such that $c\geq |v(z)|+
|u_j(w)|+1$ for all $z\in\Omega$, $w\in\partial\Omega$ and $j$. From Lemma 1
and the assumption $(dd^cv)^n\geq\mu$, it turns out that for each $j$ 
$$\int\limits_{\{u_j<v-c\}}(1-{|z|^2\over A^2})\,d\mu_j\geq 
\int\limits_{\{u_j<v-c\}}(1-{|z|^2\over A^2})\,d\mu+{1\over (n!)^2A^{2n}}
\int\limits_{\{u_j<v-c\}}(v-c-u_j)^n\,(dd^c|z|^2)^n.$$
Let $j\to\infty$, and since $||\mu_j-\mu||_\Omega\to 0$ then by Fatou Lemma
we have
$$0\geq
\liminf\limits_{j\to\infty}\int\limits_{\{u_j<v-c\}}(v-c-u_j)^n\,(dd^c|z|^2)
^n\geq
\int\limits_\Omega \liminf_{j\to\infty} \Bigl(\chi_{\{u_j<v-c\}}\,(v-c-u_j)^n\Bigr)\,(dd^c|z|^2)^n$$
$$\geq \int\limits_\Omega \chi_{\{\,\limsup
u_j<v-c\,\}}\,\bigl(\,\liminf_{j\to\infty}
|v-c-u_j|\,\bigr)^n\,(dd^c|z|^2)^n  $$  $$ \geq \int\limits_{\{\,\limsup
u_j<v-c\,\}}\,
(v-c-\limsup_{j\to\infty} u_j)^n\,(dd^c|z|^2)^n,$$
which implies that $\limsup u_j\geq v-c$ a.e. in $\Omega$ with
respect to the Lebesgue measure, and hence $\limsup u_j$ is not identically
$-\infty$ on any 
component of $\Omega$. Therefore by Corollary 7.3 in [B-T2] we can 
find a non-negative psh function $g$ in $\Omega$ such that the set
$\{g\not=\limsup u_j\}$ is pluripolar;
i.e., a set of $C_n$-capacity zero. So $g\geq v-c$ a.e. in $\Omega$ with
respect to the Lebesgue measure and hence $g$ is a bounded function in $\Omega$.
We shall prove that $u_j\to g$ in $C_n$-capacity on each compact subset $E$
of $\Omega$. For
each $\delta>0$ we have $$C_n\bigl(E\cap \{ |g-u_j|\geq\delta\}\bigr)\leq
C_n\Bigl(E\cap \bigl\{ |g-\sup\{u_j,u_{j+1},\dots\}
|\geq{\delta\over 2}\bigr\}\Bigr) $$  $$+C_n\bigl\{ |\sup\{u_j,u_{j+1},\dots\}
-u_j|\geq{\delta\over 2}\bigr\}.$$
By Proposition 5.1 in [B-T2] we know that outside a set of $C_n$-capacity
zero, the functions
$\sup\{u_j,
u_{j+1},\dots\}=\sup^\ast\{u_j,u_{j+1},\dots\}$ decrease to the function
$\limsup u_j=g$ when $j\nearrow\infty$.
Hence combining with the quasi-continuity of psh functions, Dini's theorem
implies that $\sup\{u_j,
u_{j+1},\dots\}\to g$ uniformly on $E$ outside a set of the $C_n$-capacity
less than any given constant. Thus
the first term on the right hand side of the last inequality converges to
$0$ as $ j\to\infty.$ To see 
that the second term also converges to $0$ as $ j\to\infty,$ we first show 
the following inclusion
$$\bigl\{\, |\sup\{u_j,u_{j+1},\dots\} -u_j|\geq{\delta\over 2}\,\bigr\}\subset
\bigcup_{l=0}^\infty \bigl\{\, |u_{l+j+1}-u_{l+j}|\geq{\delta\over
{2^{l+j+2}}}\,\bigr\}.$$ 
For this let $z_0$ be a point in the set on the left hand side. We choose
an integer  $l_0$ such that
$|u_{l_0+j+1}(z_0)-u_j(z_0)|\geq \delta/4$. Assume $z_0\not\in 
\cup_{l=0}^{l_0-1} \bigl\{ |u_{l+j+1}-u_{l+j}|\geq{\delta/
{2^{l+j+2}}}\bigr\}.$ Then
$$|u_{l_0+j+1}(z_0)-u_{l_0+j}(z_0)|\geq |u_{l_0+j+1}(z_0)-u_j(z_0)|-
\sum\limits_{l=0}^{l_0-1}|u_{l+j+1}(z_0)-u_{l+j}(z_0)|$$   $$\geq
{\delta\over 4} -
\sum\limits_{l=0}^{l_0-1} {\delta\over {2^{l+j+2}}}\geq {\delta\over
{2^{l_0+j+2}}}, $$
which implies $z_0\in \bigl\{ |u_{l_0+j+1}-u_{l_0+j}|\geq{\delta/
{2^{l_0+j+2}}}\bigr\}$, and hence the above
inclusion holds. So we have 
$$C_n\bigl\{ \,|\sup\{u_j,u_{j+1},\dots\} -u_j|\geq{\delta\over
2}\,\bigr\}\leq \sum\limits_{l=0}^\infty
C_n \bigl\{\, |u_{l+j+1}-u_{l+j}|\geq{\delta\over {2^{l+j+2}}}\,\bigr\}$$
$$\leq \sum\limits_{l=0}^\infty {{2^{n(l+j+2)}}\over {\delta^n 2^{(n+1)(l+j)}}}=
{4^n\over{\delta^n 2^j}}\longrightarrow 0,\qquad as\ j\to\infty.$$
Therefore $u_j\to g$ in $C_n$-capacity on each compact subset $E $ of
$\Omega$ and it then turns out from Theorem 1 that 
$(dd^cu_j)^n\to(dd^cg)^n$ as currents. Hence $\mu=(dd^cg)^n$ in $\Omega$
and the proof is complete.
\bigskip 
As a consequence of Theorem 4 we also get the following result in [C-S].
\bigskip
\noindent{\bf Corollary.} \it Suppose that $\Omega$ is a bounded domain in
$\C^n$. If there exists a function $v\in PSH(\Omega)\cap 
L^\infty(\Omega)$ such that $(dd^cv)^n\geq fd\lambda$, where the function
$f\in L^1(\Omega)$, then 
there exists a function $u\in PSH(\Omega)\cap L^\infty(\Omega)$ such that
$(dd^cu)^n=fd\lambda$.\rm
\bigskip
\noindent{\it Proof.} Since every non-negative integrable function is, in
the $L^1$-norm,
the limit of some sequence of non-negative continuous functions with
compact support, then Corollary follows from
Theorem D in [B-T1] and Theorem 4.
\bigskip
Note that without the assumption of the existence
of a subsolution, neither Theorem 4 nor Corollary is true, as can be seen
from the fact that there exists a positive measure $\mu=fd\lambda $
with $f\in L^1(\Omega)$ which is not the Monge-Amp\`ere measure of a
bounded psh function, see [C-S]. 

\bigskip \bigskip
\bigskip \noindent{\bf References } 
\bigskip\smallskip 

\noindent [B-T1] $\,$E.Bedford and B.A.Taylor, {\it The Dirichlet problem
for the complex
Monge-Amp\`ere} 

\quad {\it operator}. Invent. Math. {\bf 37} (1976), 1-44.

\noindent [B-T2] $\,$E.Bedford and B.A.Taylor, {\it A new capacity for
plurisubharmonic
functions}. Acta 

\quad Math., {\bf 149} (1982), 1-40.

\noindent [C] \quad\enskip U.Cegrell, {\it Discontinuit\'e de l'op\'erateur
de Monge-Amp\`ere complexe}.
C. R. Acad. 

\quad Sci. Paris Ser. I Math., {\bf 296} (1983), 869-871.

\noindent [C-S] $\,\ $ U.Cegrell and A.Sadullaev, {\it Approximation of
plurisubharmonic functions
and the 

\quad Dirichlet problem for the complex Monge-Amp\`ere operator}. Math.
Scand. {\bf 71}  

\quad (1993), 62-68.

\noindent [L1] \quad P.Lelong, {\it Fonctions plurisousharmoniques et
formes differentielles
positives}. Gor-

\quad don and Breach, Paris, 1968.

\noindent [L2] \quad P.Lelong, {\it Discontinuit\'e et annulation de
l'op\'erateur de Monge-Amp\`ere
complexe}.  

\quad Lecture Notes in Math., Springer-Verlag, Berlin {\bf 1028} (1983),
219-224. 

\noindent [R] \quad\enskip L.I.Ronkin, {\it Weak convergence of the current 
$[dd^cu_t]^q$ and asymptotics of the finite}

\quad {\it order function
for holomorphic mappings of regular growth}. Sibirski\u\i\enskip Mat. Zh., 

\quad {\bf 25}:4
(1984), 167-173.

\noindent [X] \enskip$\ \,$ Y.Xing, {\it On convergence of the current
$(dd^cu_t)^q$.}
To appear in S\'eminaire P.Lelong-

\quad P.Dolbeault-H.Skoda, Lecture Notes in Math.

\bigskip 
\noindent Department of Mathematics, University of Ume\aa, S-901 87 Ume\aa,
Sweden
\smallskip
\noindent Electronic mail:\enskip Yang.Xing@mathdept.umu.se

\bye